\documentclass[11pt]{article}
\usepackage{amsmath}
\usepackage{amsthm}
\usepackage{amsfonts}
\usepackage{amssymb}
\newcommand{\be}{\begin{equation}}
\newcommand{\ee}{\end{equation}}
\newcommand{\bea}{\begin{eqnarray}}
\newcommand{\eea}{\end{eqnarray}}
\newcommand{\bean}{\begin{eqnarray*}}
\newcommand{\eean}{\end{eqnarray*}}
\newcommand{\brray}{\begin{array}}
\newcommand{\erray}{\end{array}}
\newcommand{\ben}{\begin{equation}{nonumber}}
\newcommand{\een}{\end{equation}{nonumber}}

\newtheorem{dfn}{Definition}[section]
\newtheorem{thm}[dfn]{Theorem}
\newtheorem{lema}[dfn]{Lemma}
\newtheorem{pro}[dfn]{Proposition}
\newtheorem{coro}[dfn]{Corollary}
\newtheorem{xmpl}[dfn]{Example}
\newtheorem{rmrk}[dfn]{Remark}

\newcommand{\bdfn}{\begin{dfn}}
\newcommand{\bthm}{\begin{thm}}
\newcommand{\blema}{\begin{lema}}
\newcommand{\bpro}{\begin{pro}}
\newcommand{\bcoro}{\begin{coro}}
\newcommand{\bxmpl}{\begin{xmpl}}
\newcommand{\brmrk}{\begin{rmrk}}

\newcommand{\edfn}{\end{dfn}}
\newcommand{\ethm}{\end{thm}}
\newcommand{\elema}{\end{lema}}
\newcommand{\epro}{\end{pro}}
\newcommand{\ecoro}{\end{coro}}
\newcommand{\exmpl}{\end{xmpl}}
\newcommand{\ermrk}{\end{rmrk}}

\newcommand{\A}{\mathcal A}

\newcommand{\C}{\mathcal C}

\numberwithin{equation}{section}

\newcommand{\MBBZ}{\mathbb{Z}}
\newcommand{\MBBR}{\mathbb{R}}
\newcommand{\MBBT}{\mathbb{T}}
\newcommand{\MBBC}{\mathbb{C}}
\begin{document}
\begin{center}
{\large {\bf  Invariants for Normal Completely Positive  Maps on the Hyperfinite $II_1$ Factor  }}\\

{\large Debashish Goswami, Lingaraj Sahu {\footnote {The  author would like to
acknowledge the support of National Board of Higher Mathematics, DAE, India.}}}\\

\emph{Stat-Math Unit, Indian Statistical Institute,\\ 203,
B.T. Road, Kolkata-700108, India.
\\ email
:  goswamid@isical.ac.in,
lingaraj@gmail.com }\\
\vspace{1cm}
 \bf{ Dedicated to Prof K.B. Sinha}
\end{center}
\begin{abstract}
We investigate  certain classes of  normal completely positive (CP) maps  on  the 
hyperfinite $II_1$ factor $\mathcal A$. Using the representation theory of a suitable irrational rotation algebra, we propose some computable invariants for such CP maps.  
\end{abstract}
\section{Introduction}
Quantum dynamical semigroups (QDS), i.e. $C_0$-semigroups of completely positive (CP)  contractive  maps on $C^*$ or von Neumann algebras (assumed to be normal in the von Neumann algebra case), are important objects both from physical and mathematical viewpoint. It is thus quite natural to look for a nice classification scheme for them. As in every branch of mathematics, one would like to find one or more computable (preferably numerical) invariants for QDS with respect to some suitable equivalence relation, which is by now accepted as the so-called cocycle congugacy introduced and studied in a series of papers by   Arveson (\cite{Arv1},\cite{Arv2},\cite{Arv3}), Powers \cite{Pow}, Bhat (\cite{Bh1},\cite{Bh2}) and others.
There has already been a considerable amount of literature on this problem as well as the related (and in some sense equivalent) problem of classifying product system of Hilbert modules, thanks to the intesive works by a number of mathematicians including  Arveson, Bhat, Skeide, Tsirelson, Barreto, Liebscher, Verschik, to name only a few (see, e.g.  \cite{Arv1}, \cite{Sk_etal}, \cite{Sk2} and references therein). However, while much is known about QDS on ${\mathcal B}(\mathcal H)$ (where $\mathcal H$ is a separable Hilbert space), and a nice numerical invariant (`Arveson index', see, e.g.  \cite{Arv2} )is available in this case, QDS on other types of von Neumann algebra are not so well understood, in the sense that there is  not yet any satisfactory numerical (or  easily computable) invariant for QDS (or, equivalently, product systems of Hilbert modules) on general von Neumann algebras. In order to construct such invariants, it is reasonable to first  restrict attention  to just a single CP map insetead of semigroup of CP maps. 
 This is what we attempt to do in the present article for the hyperfinite type $II_1$ factor.  While we are not yet able to come up with a satisfactory numerical invariant for an arbitray normal CP map on the hyperfinite $II_1$ factor, we do get a nice invariant for an interesting class of CP maps, namely the `pure' ones. This invariant, though not numerical, is given by a quadruplet $(m,\nu, b_1,b_2)$ where $m$ is a nonnegative integer (so can be thought of as an analogue of `Arveson index' in this case), $\nu$ is a measure (with suitable property) on the two-torus $ \MBBT^2$, and $b_1,b_2 :  \MBBT^2 \rightarrow M_m ( \MBBC)$ are unitary valued $\nu$-measurable maps. Thus, the invariant is in some sense not a very abstract object, and easily computable. It can  be described at least by a sequence of complex numbers, if we describe the measure $\nu$ as well as the bounded functions $z \mapsto b_i(k,l):=(k,l){ \rm -th~entry~of~}~b_i(z),~i=1,2,$ in terms of the  Fourier coefficients. 
 Apart from the pure CP maps, we also study in some details another class which we call `extendible'. We give necessary and sufficient conditions for extendibility, and propose some invariants for such CP maps, including numerical invariants.   

\section{Preliminaries}

 \subsection{Intertwiner Module}
Let $\mathbf h$ be a Hilbert space, $\A \subseteq \mathcal B( \mathbf h)$ be a von Neumann algebra. Given  a normal unital CP map $ T :\mathcal A\rightarrow {\mathcal A} $, by Stinespring's and KSGNS theorems (see \cite{Lan})    we can obtain  a Hilbert von Neumann  $\mathcal A-\mathcal A$ 
bimodule (to be called KSGNS bimodule) $\mathcal{E}_T\subseteq \mathcal B(\mathbf h,{\mathcal K}_T)$  for some Hilbert space ${\mathcal K}_T,$ a normal  representation  $\pi_T:\mathcal A \rightarrow {\mathcal {L}}({\mathcal{E}_T})\subseteq \mathcal B(\mathcal K)$
and $\xi_T \in \mathcal{E}_T$ such that 
 $T(a)=\langle \xi_T,
\pi_T(a)\xi_T \rangle,\forall a \in \mathcal A.$  Here $\langle \xi,
\eta \rangle:=\xi^* \eta$ is the  $\A$-valued inner product on $\mathcal{E}_T$ and $\pi_T(\A) \xi_T \A$ is total in $\mathcal{E}_T$ which implies that 
$\pi_T (\A)\xi_T \mathbf h$ is total in $ {\mathcal K}_T.$ 
\noindent Setting $\rho_{T}:\mathcal A^\prime
\rightarrow\mathcal B({\mathcal K}_T)$ by
$$ \rho_T(a^\prime) \pi_T (a)\xi h= \pi_T (a)\xi_T a^\prime h,$$ we get a normal representation $\rho_T$  of ${\mathcal A}$ such that $\pi_T(\mathcal A) $
and $\rho_T(\mathcal A^\prime~)$ commute and the Hilbert (right) module ${\mathcal E}_T$ coincides with   
the intertwiner module $I(\rho_T,\mathcal A^\prime~)$ (c.f. \cite{Gos},\cite{Ski}) of $(\rho_T,\mathcal A^\prime~)$, which is given by :  
\[I(\rho_T,\mathcal A^\prime~)= \{ R\in \mathcal B(\mathbf h,{\mathcal K}_T):{R}a^\prime~=
\rho_T(a^\prime)~R,
 \forall a^\prime\in {\mathcal A}^\prime \}.\]

 \blema Let $T$ and ${\mathcal E}_T$ be discussed as above. Then
 \[{\mathcal {L}}({\mathcal{E}}_T)=(\rho_T(\A^\prime))^\prime\]
 \elema
 \begin{proof}
Let  $B\in \mathcal{L}(\mathcal E_T).$ For $a\in \A, a^\prime \in \A^\prime, h\in \mathbf h,$ we have 
\[B \rho_T(a^\prime) \pi_T(a) \xi_T h=B \pi_T(a)~\xi_T a^\prime h=  \rho_T(a^\prime)(  B \pi_T(a)\xi_T) h.\]
Since $\pi(\A)\xi_T \mathbf h$ is total in ${\mathcal K}_T,$ 
$B \in (\rho_T(\A^\prime))^\prime.$

\noindent For the converse, let $B \in (\rho_T(\A^\prime))^\prime.$  For $ a^\prime \in \A^\prime$ and $\xi\in \mathcal E_T$ we have 
\[B \xi  a^\prime  =B \rho_T(a^\prime) \xi =\rho_T(a^\prime) B \xi.\]
This implies that $B \xi$ is in the intertwiner module $\mathcal E_T.$
\end{proof}

 Let us now introduce various equivalence relations on the set CP($\mathcal A)$ of normal unital CP maps. From now on, unless there is a specific need to mention, we shall omit the adjectives `normal CP unital'. 
 
 The following equivalence relation is motivated by the definition of cocycle conjugacy for QDS.
 \bdfn 
  We say that two elements $T,S$ of CP($\mathcal A$) are equivalent, denoted as $T \approx S$, if there exist  unitaries $u,v$  in $\mathcal A,$  such that for all $a \in {\mathcal A}$, 
$$ S(a)=v^*T(u^* au)v.$$
\edfn

We also introduce another (stronger) equivalence relation $\sim$.
\bdfn
$T \sim S$ if there exists unitary $u \in {\mathcal A}$ such that $S(a)=T(u^*a u)$ for all $a \in {\mathcal A}$.
\edfn

If $S \approx T$ as above,  we have
\bean
&&S(a)=v^*T(u^* au)v=v^* \langle \xi_T,
\pi_T(u^* au)\xi_T \rangle v\\
&&=
v^*  \xi_T^*
\pi_T(u)^*  \pi_T(a) \pi_T(u)\xi_T v
=\langle \xi_S,
\pi_T(a)\xi_S\rangle,
\eean
where $\xi_S= \pi_T(u)\xi_T v .$ So the corresponding KSGNS bimodules can be chosen to be the  same,
say, ${\mathcal E}$, and   $\xi_S=W \xi_T$ for some unitary $W\in {\mathcal L}(\mathcal E).$

It is not unreasonable to identify two CP maps if the KSGNS bimodules are isomorphic, i.e. can be chosen to be the same. This leads to the third equivalence relation.
\bdfn
We say that $T$ and $S$ in CP($\mathcal A$) are KSGNS-equivalent, to be denoted by $T \stackrel{\prime}{\sim} S$, if the bimodules ${\mathcal E}_T$ and ${\mathcal E}_S$ are isomorphic.
\edfn
It is clear from our discussion that $T \sim S \Rightarrow T \approx S \Rightarrow T \stackrel{\prime}{\sim} S$. However, the converse implications can be shown to be false in general.

\subsection{Second order irrational rotation algebra and its representations}
Fix an irrational number $\theta$, and consider the irrational rotation algebra   ${\mathcal A}_{\theta},$ (c.f. \cite{Dav}), which is 
the universal $C^*$-algebra
generated by two  unitaries  $U$  and  $V $ satisfying the   Weyl commutation relation
$UV=\lambda VU,$ where $\lambda=e^{2\pi i \theta}$.  Let $ \tau $ be 
the unique faithful normalized trace on $\mathcal A_{\theta}$
and  $ (\mathcal H _ {\tau} \equiv L^2(\mathcal A_\theta, \tau) ,\pi_{\tau}, $ 1 $),$ be the associated GNS triple.  From now, we fix the von Neumann algebra  $ \mathcal A $ as  the (unique upto isomorphism) hyperfinite type  $II_1 $ factor,  identified with the weak closure of   $\mathcal A _
 \theta$ in $ \mathcal B(\mathcal H_\tau).$ It is well-known that $\A_\theta$ is isomorphic with the crossed product ${\mathcal C}(\MBBT) \rtimes \MBBZ$, with respect to the action   described in \cite{Dav}. Define the unitary operators
 $U_l,V_l,U_r$ and $ V_r$ in $ \mathcal B(\mathcal H_\tau)$ by
$U_l a=Ua,~V_la=Va,~U_ra=aU$ and $V_ra=aV$ for  $ a\in \mathcal A_{\theta} \subseteq {\mathcal H}_\tau.$ It is clear that the $C^*$-algebras  $C^*(U_l,V_l)$ and $
C^*(U_r,V_r)$ are weakly dense $\ast$-sub-algebras of $\A$ and $\A^\prime$
respectively. The vector  $1$ is cyclic and separating  both for  $\A$ and $\mathcal
A^\prime.$
 
\noindent Consider the  $\ast$- sub-algebras  $\A^{\rm l,fin}$ and  $\A^{\rm r,fin}$ of  $\A$  and $\A^\prime$ 
respectively, where 
$ \A^{\rm l,fin}$ is the unital $\ast$-algebra generated by all  polynomials in the unitaries $ U_l$ and $V_l$, and  
$\A^{\rm r,fin}$ is the similar algebra obtained by replacing $U_l,V_l$ by $U_r,V_r$ respectively. Clearly, $\A^{\rm l,fin}$  is norm-dense in $\A_\theta$. Moreover, let $\A^{\rm fin}$ denote  the  $\ast$-algebra of $\A^{\rm l,fin} \otimes_{\rm alg} \A^{\rm r,fin}$ generated by polynomials in all the four unitaries $U_l, V_l, U_r$ and $V_r$ and $\mathcal B$ be the $C^*$ algebra obtained by completing $\A^{\rm fin}$ in the norm-topology of ${\mathcal B}({\mathcal H}_\tau \otimes {\mathcal H}_\tau)$.   
\blema \label{2theta}
The $C^*$-algebra  $\mathcal B$ is  isomorphic (as $C^*$ algebra) with the crossed product $C^*$-algebra   ${\C}(\MBBT^2)\rtimes \MBBZ^2$  with respect to  the  action $\alpha$ given by
$\alpha(k_1,k_2)(f)(z_1,z_2)=f(\lambda^{-k_1} z_1,\lambda^{-k_2}z_2).$\\
\elema
\begin{proof}
By definition,  the crossed product $C^*$-algebra   ${\C}(\MBBT^2)\rtimes \MBBZ^2$ with respect to  the  action $\alpha$ is the universal $C*$ algebra generated by  four unitaries $W_1,W_2,W_3,W_4$ satisfying $$ W_1 W_2=W_2 W_1,~W_3W_4=W_4W_3,~ W_1W_3=\lambda W_3W_1,~ W_2W_4=\lambda W_4 W_2.$$  
It can be vriefied by a straightforward calculation that  the linear map  $\eta:{\C}(\MBBT^2)\rtimes \MBBZ^2\rightarrow \A^{\rm fin}$ defined by \\
$$ \eta (W_1)=U_l \otimes 1,~
\eta (W_2)=1 \otimes V_r,~
\eta (W_3)=V_l^{-1}\otimes V_r,~
\eta (W_4)=U_l^{-1}\otimes U_r,$$
gives   the required isomorphism.
\end{proof}
From the general theory of $C^*$ crossed product, it is easy to show that $\mathcal B \cong \C(\MBBT^2) \rtimes \MBBZ^2$ is simple, and thus any nonzero representation is a $C^*$ isomorphism. Thus, the $C^*$ algebras $\mathcal B$, $\C(\MBBT^2) \rtimes \MBBZ^2$ and $C^*(U_l,V_l,U_r,V_r)(\subseteq {\mathcal B}({\mathcal H}_\tau))$ are all canonically isomorphic, and we shall identify them whenever there is no chance of confusion, and shall call this $C^*$ algebra the `second order irrational rotation algebra', and denote it by  $\mathcal A_\theta^{(2)}.$ It is actually the $4$-dimensional noncommutative torus   in the sense of Rieffel and Schwarz \cite{Rief}, corresponding to the skew-symmetric $4 \times 4$ matrix  $A=(( a_{ij} ))$, with $a_{12}=a_{34}=\theta, a_{21}=a_{43}=-\theta$, and all other entries are zero. Moreover, from the proof of the above Lemma, it is clear that the map $\eta$ gives an algebraic isomorphism between $\A^{\rm fin}$ and the $\ast$-algebra generated by 
 $W_i$, $i=1,...,4.$

 The weak closure of the $C^*$-algebra $\mathcal A_\theta^{(2)}$  in the weak topology inherited from ${\mathcal B}({\mathcal H}_\tau \otimes {\mathcal H}_\tau)$ is   $\mathcal A \otimes \mathcal A^\prime \cong
W^*(U_l,V_r)\rtimes \MBBZ^2 \cong L^\infty(\MBBT^2) \rtimes \MBBZ^2 ,$ which is again a concrete realization of  the (unique upto isomorphism) $II_1$ hyperfinite factor,
with the faithful normal trace $\tau^{(2)}$  given by the vector state  $\langle 1 \otimes 1, ~\cdot~ 1 \otimes 1 \rangle$. We introduce here some notation for future use. We denote by $\xi_0$ the vector $1 \otimes 1$, and by $c^{(r)}$ (where $c \in \A \otimes \A^\prime$) the  operator of right multiplication by $c$ on ${\mathcal H}_\tau \otimes {\mathcal H}_\tau$, i.e. $c^{(r)}b:=bc$, for $b \in \A \otimes \A^\prime$, viewed as an $L^2$-element. The commutant of $\A \otimes \A^\prime$, i.e. $\A^\prime \otimes \A$, coincides with the set $\{ c^{(r)}~:~c \in \A \otimes \A^\prime \}$. We shall use the notation $\tau^{(2)}(X)$ for $\langle \xi_0, X \xi_0 \rangle$ whenever it makes sense, even when the operator $X$ is unbounded but contains $\xi_0$ in its domain.

  Now we make an important observation  regarding extension of representations.
 \bpro \label{poly-cont}
Any representation  $\pi$ of the $*$-algebra $\A^{\rm l,fin}$  ($\A^{\rm fin}$) extends to the $C^*$-algebra $\mathcal A_\theta$  ($\mathcal A_\theta^{(2)}).$
 \epro
 \begin{proof}
 Let us give a sketch of the proof for $\A^{\rm l,fin}$ only. The proof for $\A^{\rm fin}$ is similar.
 Let $\pi : \A^{\rm l,fin} \rightarrow {\mathcal B}({\mathcal H^\prime})$ be a representation. Thus, we have two unitaries $\pi(U_l)$ and $\pi(V_l)$ acting on ${\mathcal B}({\mathcal H^\prime})$, satisfying the commutation relation $\pi(U_l)\pi(V_l)=\lambda \pi(V_l)\pi(U_l)$.  
 Since the commutative $C^*$-algebra $\mathcal C(\MBBT)$  is the universal $C^*$-algebra generated  by a unitary, we can get a representation $\pi_0$, say, of $\mathcal C(\MBBT)$ which maps $f \in \mathcal C(\MBBT)$ to $f(\pi(U_l))$. Clearly, $(\pi_0,\pi(V_l))$ is a covariant representation for the action $(m \cdot f)(z):=f(\lambda^m z)$, and thus, there exists a representation $\tilde{\pi} : \mathcal C(\MBBT) \rtimes \MBBZ $ ($\equiv {\mathcal A}_\theta$) into $\mathcal B(\mathcal H^\prime)$ satisfying $\tilde{\pi}(U_l)=\pi(U_l),\tilde{\pi}(V_l)=\pi(V_l)$.  
 \end{proof}

 \subsection{ Irreducible representations of  ${\mathcal A}_\theta^{(2)}$}
 Following \cite{Bre}, we say that a separable representation $\pi$ of the irrational rotation algebra
$\mathcal A_\theta$ has the  uniform multiplicity $m$
 if the restriction of $\pi$ to $\mathcal
C(\MBBT^2)$ has uniform multiplicity $m\in \{1,2,\cdots \infty\}.$
The factor representations (in particular,  the irreducible
representations) are in above class.
 In \cite{Bre}, the author  gave a nice 
invariant of uniform multiplicity representation on the irrational
rotation algebra $\A_\theta$. These representations  are
classified by their multiplicities $m$,
 a quasi-invariant Borel measure $\nu$ (where quasi-invariance is with respect to rotation by the angle
  $2\pi \theta$) and  a unitary one-cocycle   $b$. It is not difficult to verify that the arguments of \cite{Bre} 
     can 
  easily be extended  to  the $2$-nd order irrational rotation algebra ${\mathcal A}_\theta^{(2)}$. We summarize this fact as a theorem below, the proof of which is omitted since it is very similar to the proofs of analogous results in \cite{Bre}.
\bthm
\label{invcp}
(i) Any 
    irreducible representation of ${\mathcal A}_{\theta}^{(2)}=C(\MBBT^{2}) \rtimes \MBBZ^2=C^{*}(W_{1},
W_{2},W_{3},W_4)$ is unitarily equivalent to a representation of the form   $\pi_{(m,\nu,b_1,b_2)}$ described below, where   $m \in \mathbb{N} \cup \{0\}$,  $\nu$ is an ergodic regular Borel measure  on $\MBBT^2$ which is also quasi-invariant with respect to the action $\alpha$ of $\MBBZ^2$ given by $\alpha_{(n_{1},n_{2})}z=(\lambda^{n_{1}}z_{1},\lambda^{n_{2}}z_{2})$, and two $\nu$-measurable unitary valued functions $b_{1},b_{2}:\MBBT^2 \rightarrow \mathcal{M}_{m}(\mathbb{C})$ such that $\pi_{(m,\nu,b_1,b_2)} : {\mathcal A}_\theta^{(2)} \rightarrow {\mathcal B}(L^2(\MBBT^2,\nu)\otimes \MBBC^m) $ is given by:
\begin{equation}
\begin{split}
\pi(W_{1})f(z)=&z_{1}f(z)\\
\pi(W_{2})f(z)=&z_{2}f(z)\\
\pi(W_{3})f(z)=&b_{1}(z)\sqrt{\frac{d\nu_{1,0}}{d\nu}}(z)f(\lambda z_{1},z_{2})\\
\pi(W_{4})f(z)=&b_{2}(z)\sqrt{\frac{d\nu_{0,1}}{d\nu}}(z)f(z_{1},\lambda z_{2}),
\end{split}
\end{equation}
where $f : \MBBT^2 \rightarrow M_m(\MBBC)$, and $ \frac{d \nu_{n_{1},n_{2}}}{d\nu}$ is the Radon-Nikodym derivative of the measure $\nu_{n_{1},n_{2}}$ given by 
\[
\nu_{n_{1},n_{2}}(E)=\nu \{\alpha_{n_{1}n_{2}}z:z \in E\}
\]
with respect to $\nu$. Note that $\nu_{n_{1},n_{2}}$ is equivalent   to  $\nu$ as $\nu$ is quasi-invariant.\\
(ii)
Two irreducible representations  $\pi_{(m,\nu,b_{1},b_{2})}$ and $\pi_{(\tilde{m},\tilde{\nu},\tilde{b_{1}},\tilde{b_{2}})}$. are unitarily equivalent  if and only if
\begin{enumerate}
\item $m=\tilde{m}$ and $\nu=\tilde{\nu}$ can be chosen;\\
\item there  exists a  unitary-valued, $\nu$-measurable map $W:\MBBT^2 \rightarrow \mathcal{M}_{m}(\mathbb{C})$ satisfying $$ W(z)b_{1}(z)=\tilde{b_{1}}(z)W(\lambda z_{1},z_{2}),~~
W(z)b_{2}(z)=\tilde{b_{2}}(z)W(z_{1},\lambda z_{2})
.$$
\end{enumerate}
\end{thm}

\section{Invariants for CP maps on type $II_1$ factor}

\subsection{State associated with CP map}
Given a normal CP unital map  $T$ on the type $II_1$ factor $\mathcal A,$ consider the associated KSGNS bimodule ${\mathcal E}_T$, cyclic element $\xi_T \in {\mathcal E}_T$, and the representations $\pi_T,\rho_T$ discussed in  section 2.  
Now we define  a representation  $\tilde{\pi}_T:{\mathcal A } \otimes _ {alg}
\mathcal A^\prime \   {\rightarrow}{\mathcal B}(\mathcal K_T)$ by
$$\tilde{\pi}_T(a\otimes a^\prime~)= \pi_T (a)\rho_T(a^\prime~).$$
This is a homomorphism since $\pi_T(a)$ and $\rho_T(a^\prime)$ commute.  
By Proposition \ref{poly-cont}, 
the restriction of   $\tilde{\pi}$ to the weakly dense subalgebra ${\mathcal A}^{\rm fin}$   extends to the  $C^*$-algebra $\A_\theta^{(2)},$ and we denote this extension again by the same notation.  
Since $\pi_T({\mathcal A})\xi_T{\mathcal H}_\tau$  is dense in ${\mathcal K}_T,$ and $1$ is cyclic for $\mathcal A^\prime$, it is clear that $\xi_T 1 \in $ is cyclic
for ($\tilde{\pi}_T,{\mathcal A } \otimes_{alg} {\mathcal A^\prime~},{\mathcal K}_T$). 
 We define  a  state $ \psi_T $ on ${\mathcal A}_\theta^{(2)}$ by setting 
 \be
 \label{psidef}
 \psi_T( \cdot)=\langle \xi_T 1, \tilde{\pi}_T(\cdot) \xi 1 \rangle.\ee
  Since $\tilde{\pi}_T$ is defined also on $ {\mathcal A } \otimes _{alg}
 \mathcal A^\prime~$, we shall define $\psi_T$ on  $ {\mathcal A } \otimes _{alg}
 \mathcal A^\prime~$ too, given by the same expression as in (\ref{psidef}).
 Note that we have, $\psi_T (a\otimes a^\prime~)
 =\left\langle 1,T(a)a^\prime ~1\right\rangle$. It is clear from the construction of $\psi_T$ that $({\mathcal K}_T, \tilde{\pi}_T, \xi_T~1)$ is a choice of GNS triple for $\psi_T$. Thus, for two CP maps $T,S \in {\rm CP}({\mathcal A})$, the corresponding state $\psi_T$ and $\psi_S$ have unitarily equivalent GNS representations if and only if $\tilde{\pi}_T$ and $\tilde{\pi}_S$ are unitarily equivalent. This leads to the following observation.
 \bpro
 \label{123}
 $T \stackrel{\prime}{\sim} S$ if and only if $\psi_T$ and $\psi_S$ have unitarily equivalent GNS representations, i.e. $\tilde{\pi}_T$ and $\tilde{\pi}_S$ are unitarily equivalent. 
 \epro
 \begin{proof}
 Suppose first that $T \stackrel{\prime}{\sim} S$, and ${\mathcal U}$ is an ${\mathcal A}$-linear unitary from ${\mathcal E}_T$ to ${\mathcal E}_S$ satisfying ${\mathcal U} \pi_T(a) {\mathcal U}^*=\pi_S(a)$ for all $a \in {\mathcal A}$. Using the fact that $\{ ev~:~e \in {\mathcal E}_T,~v \in {\mathcal H}_\tau\}$ and  $\{ ev~:~e \in {\mathcal E}_S,~v \in {\mathcal H}_\tau\}$ are total in ${\mathcal K}_T$ and ${\mathcal K}_S$ respectively, it is easy to verify that the map $U : {\mathcal K}_T \rightarrow {\mathcal K}_S$ defined by $$ U(ev):={\mathcal U}(e)v,~~e \in {\mathcal E}_T,~v \in {\mathcal H}_\tau,$$ extends to a unitary in the Hilbert space sense, and $U \pi_T(a) U^*=\pi_S(a),$ $U \rho_T(a^\prime)U^*=\rho_S(a^\prime)$ for all $a \in {\mathcal A},a^\prime \in {\mathcal A}^\prime$, which means that $U \tilde{\pi}_T(\cdot)U^*=\tilde{\pi}_S(\cdot)$.
 
 Conversely, given a unitary $U : {\mathcal K}_T \rightarrow {\mathcal K}_S$ such that $U \tilde{\pi}_T(\cdot)U^*=\tilde{\pi}_S(\cdot)$, we can define $\mathcal U : {\mathcal E}_T \rightarrow {\mathcal E}_S$ by $$ {\mathcal U}(e)=Ue,$$ where $e \in {\mathcal E}_T $ is viewed as an element of ${\mathcal B}({\mathcal H}_\tau, {\mathcal K}_T)$. This is easy to prove that $\mathcal U$ is an ${\mathcal A}$-linear unitary map, which also intertwins the left actions on the two modules, and thus ${\mathcal E}_T$ and ${\mathcal E}_S$ are isomorphic as bimodules. 
 \end{proof}

 We now want to chracterize the states on ${\mathcal A}_\theta^{(2)}$ which are of the form $\psi_T$ for some $T \in {\rm CP}(\mathcal A)$. 
 Let us denote by $ \mathcal S_\tau$  the  set of all states $\psi$ on
 $\A_\theta^{(2)}$ such that \\
 (i) $\psi (1 \otimes a^\prime)=\langle 1, a^\prime 1 \rangle \equiv \tau(a^\prime)$ for all $a^\prime \in {\mathcal A}^{\rm r,fin}$; \\
 and \\
 (ii) the restriction of the state $\psi$ on $\A^{\rm l,fin} \otimes 1$ admits a  normal extension to $\mathcal A$.
 \bpro
 The   map $T \mapsto \psi_T $  is a bijection from $ CP(\mathcal A)$ to $ \mathcal S_\tau$.
 \epro
\begin{proof}
It is clear that for $T \in {\rm CP}(\mathcal A),$ $\psi_T$ is an element of $\mathcal S_\tau$.  Let us first show that it is one-to-one.
Let $T$ and $T^\prime$ be two CP maps such that $\psi_T=\psi_{T^\prime}.$
For any $a\in \A^{\rm l,fin} $ and $a^\prime, b^\prime \in \A^{\rm r,fin}$, we have
\bean
&&\psi_T (a\otimes a^\prime~ b^\prime)\\
&& =\left\langle  1,T(a) (a^\prime)^* b^\prime 1\right\rangle\\
&& =\left\langle  1,T^\prime(a)(a^\prime)^* b^\prime 1\right\rangle.
\eean
This gives
\bean
&& \left\langle a^\prime 1,T(a) b^\prime 1\right\rangle=\left\langle a^\prime 1,T^\prime(a) b^\prime 1\right\rangle .
\eean
Since $\{a^\prime 1: a^\prime \in \A^{\rm r,fin} \}$ is total in $\mathcal {H}_\tau$
it follows that $T(a)=T^\prime(a),~\forall a\in \A^{\rm l,fin}$, hence for all $ a \in \A$. 

Now we show that the map $T \mapsto \psi_T$ is onto. 
Let $\psi  \in  \mathcal S_\tau.$  We define a sesquililear form  $T(a)$    for $a$ in the completion of $ \A^{l,fin}$, i.e. $\A_\theta$, by the following, where  $a^\prime_1,a^\prime_2 \in \A^{\rm r,fin}$:
$$ \langle a_1^\prime~1, T(a) a_2^\prime~1 \rangle= \psi(a \otimes (a_1^\prime)^* a_2).$$
\noindent For $a_i, i=1,2\cdots n$ in $\A^{\rm r,fin}$ and $b_i,i=1,2\cdots n$ in $\A^{l,fin}$, we have,
 \bea
 \label{positive}
 \sum_{i,j} \langle a_i,T(b_i^*b_j) a_j\rangle  
 = \sum_{i,j} \psi(b_i^*b_j\otimes a_i^* a_j).
  \eea
 Since   $((b_i^*b_j))$ and $(( a_i^* a_j))$ are positive in $M_n(\A)$ and $M_n(\A^\prime)$ respectively, and each $a_i$ commutes with each $b_j$, it follows by standard arguments, similar to those used for proving the positivity of Schur product of two positive (numerical) matrices, that the right hand side of (\ref{positive}) is positive.  By the condition $\psi|_{1 \otimes \A^{\rm r,fin}}=\tau$,  $T(1)=1$. This, combined with the positivity of the right hand side of (\ref{positive}), suffices to show that  $T(a)$ extends to a bounded  map
  from $\A_\theta$ to $\mathcal B(\mathcal H_\tau)$, to be denoted by the same symbol again. It is also easy to see that $T$ is CP and  $T(\A_\theta)\subseteq (\A_\theta)^{\prime \prime}=\A$ 
 Finally, as the resriction of the state $\psi$ on $\A_\theta \otimes 1$ is normal it follows that  the map $T$ extends as  a normal CP unital map on $\A$.
\end{proof}

\subsection{Invariant for a   pure CP map}
\bdfn
An element $T\in {\rm CP}(\mathcal
A)$ is said to be a pure CP map (see \cite{Tsu}, \cite{Arv}) there does not exist any CP (not necessarily unital) normal map $S$ on $\A$, other than scalar multiples of $T$, such that $T-S$ is CP.
\edfn

 \bpro A CP map $T:\A\rightarrow \A$ is  pure  if and only if the  state $\psi_T$ is  a pure state
on the  $2$-nd order  irrational rotation  $C^*$-algebra
$\A_\theta^{(2)}.$ 
\epro
 \begin{proof}
First  we note that the state $\psi_T$ is  a  pure state if and only if the GNS representation $\tilde{\pi}_T$ is irreducible. Since $\tilde{\pi}_T$ is constructed out of two commuting representations $\pi_T$ and $\rho_T,$ the irreducibility  of $\tilde{\pi}_T$ is  translated into the fact that  $\pi_T(\A)^\prime \bigcap (\rho_T(\A^\prime))^\prime=\pi_T(\A)^\prime \bigcap \mathcal L(\mathcal E_T)$
is trivial.

Now, it suffices to prove that    $ T$ is pure if and only if $\pi_T(\A)^\prime \bigcap \mathcal L(\mathcal E_T)$ is trivial. This fact can be deduced from Corollary 3.7 of \cite{Tsu} by some straghtforward argumemts, or by a direct argument along the lines of the proof of a similar fact for states as in \cite{Tak}.
 \end{proof}

Thus, 
given  any pure   CP map  $T$ on $\A,$ the  state $\psi_T$ is  a pure state
on the  $2$-nd order  irrational rotation  $C^*$-algebra
$\A_\theta^{(2)},$ i.e.   
 its GNS representation  $  \tilde{\pi}_T$ is
irreducible. Theorem \ref{invcp} provides a quadruplet $(m,\nu,b_{1},b_{2})$, which   is
an invariant for the  CP map $T$ under the KSGNS equivalence $\stackrel{\prime}{\sim}$. Moreover, it is a complete invariant for $\stackrel{\prime}{\sim}$ by Proposition \ref{123}. It is also easy to see that if $T$ is pure and $S$ is KSGNS equivalent to $T$, then $S$ must be pure. 

For a measure $\mu$ on $\MBBT^2$, let us denote by $\mu_1$ and $\mu_2$ the marginals given by: $$ \mu_1(\Delta):=\mu(\Delta \times \MBBT), ~~ \mu_2(\Delta):=\mu(\MBBT \times \Delta);$$ where $\Delta$ is a Borel subset of $\MBBT$. The resrtictions of the representations $\tilde{\pi}$ to the $C^*$-algebra 
$\A_\theta=C^*(U_l,V_l)$  and $ C^*(U_r,V_r)$ extend to normal representations of $\A$ and $\A^\prime$ (respectively) by our assumption, so it is clear that the marginals $\nu_1$  and $\nu_2$ of the  quasi invariant Borel measure $\nu$ on $\MBBT^2$ obtained above  are absolutely continuous with respect to  Lebesgue mesure on $\MBBT.$ In fact, as $\psi_T$ coincides with the canonical trace $\tau$ on $\A^\prime$, $\nu_2$ should be equivalent to the Lebesgue measure on the one-dimensional torus. 

We summarise the above discussion in form of a theorem below.
\bthm
A complete invariant for the set of pure CP maps on $\A$ under the KSGNS equivalence $\stackrel{\prime}{\sim}$ is given by the quadruplet $(m,\nu,b_1,b_2)$ as described in Theorem \ref{invcp}, with the additional conditions that $\nu_1$ is absolutely continuous with respect to the Lebesgue measure and $\nu_2$ is equivalent to the Lebesgue measure.
\ethm 

\brmrk
It should be emphasized  to avoid any confusion that the invariant discussed above  is well-defined only after fixing the isomorphism between the 
weak closure of $\mathcal A_\theta$ and the hyperfinite $II_1$ factor.  This 
is quite natural from the viewpoint of bimodule theory over factors, 
but this is different from the case of classification theory of 
automorphisms of a $II_1$ factor. 
\ermrk

\subsection{Extendible  CP maps}
 It should be noted that the state $\psi_T$ associated to a CP map $T$ may not extend to  the von Neumann algebra $\A \otimes \A^\prime$, since given a  pair of normal representations $\pi,\rho$ of $\A$ and $\A^\prime$ respectively (acting on the same Hilbert space) such  that $\pi(a)$ and $\rho(a^\prime)$ commute for all $a,a^\prime$, it is not in general possible to get a normal  representation $\phi$ of $\A \otimes \A^\prime$ such that  the restrictions of $\phi$ to $\A$ and $\A^\prime$ are respectively $\pi$ and $\rho$. In the appendix, we have given a counterexample which justifies this remark. Now, we are going to investigate when 
$\psi_{T}~$ extends to the type $II_1$ factor ${\mathcal A } \otimes  {\mathcal
A^\prime}~$ as  a normal state. Let us call the map $T$ {\bf extendible} if $\psi_T$ extends to a normal state on $\A \otimes \A^\prime$, or equivalently, $\tilde{\pi}_T$ extends to a normal representation of $\A \otimes \A^\prime$.
To give a necessary and sufficient criterion for extendibility, we need the following general result.
\bpro
\label{extend}
Let $\mathcal C \subseteq {\mathcal B}({\mathcal H_0})$ be a $C^*$ algebra, and $\tilde{\mathcal C}$ be its weak closure. Let a discrete group $\Gamma$ admit a strongly continuous unitary representation on $\mathcal H_0$ given by, say, $\Gamma \ni \gamma \mapsto u_\gamma$. Assume furthermore that $u_\gamma {\mathcal C} u_\gamma^* =\mathcal C$, so that $\beta_\gamma $ defined by $\beta_\gamma(\cdot)=u_\gamma ~\cdot ~u_\gamma^*$ defines a $\Gamma$-action on $\mathcal C$ as well as on $\tilde{\mathcal C}$. Consider the $C^*$ crossed product ${\mathcal C} \rtimes \Gamma$ and the von Neumann crossed product $\tilde{\mathcal C} \rtimes \Gamma$. Let $\phi : {\mathcal C} \rtimes \Gamma \rightarrow \MBBC$ be a state. Then, $\phi$ extends to a normal state on $\tilde{\mathcal C} \rtimes \Gamma$ if and only if for every $\gamma \in \Gamma$, the bounded linear functional $\phi_\gamma$ defined by $\phi_\gamma(c)=\phi(c \delta_\gamma)$ is weakly continuous (where $ c \delta_\gamma$ denotes the $\mathcal C$ -valued function given by $c \delta_\gamma(\gamma^\prime)=0$ if $\gamma^\prime \ne \gamma$, and is $c$ if $\gamma^\prime=\gamma$). 
\epro 
\begin{proof} 
Let us prove only the nontrivial part, i.e, the `if' part. Consider the GNS triple $(\mathcal K, \pi, \xi)$ for the state $\phi$, and denote by $\pi_0$ and $U_\gamma$ the restriction of $\pi$ to $\mathcal C$ and $\pi(\delta_\gamma)$ respectively. It is enough to prove that $\pi_0$ is normal. Consider a directed family $0 \leq c_\alpha \uparrow c $, where $c_\alpha,c \in \mathcal C$. We have to show that \be
\label{678}
\langle w, \pi_0(c_\alpha)w \rangle \uparrow \langle w, \pi_0(c) w \rangle~~\forall w \in {\mathcal K}. \ee Since the family $\pi_0(c_\alpha)$ is uniformly bounded (in operator norm) by $\| \pi_0(c) \|$ (as $0 \leq \pi_0(c_\alpha) \leq \pi_0(c)$), it is enough to verify (\ref{678}) for $w$ belonging to the dense set $\mathcal D$ spanned by vectors of the form $\pi_0(c)U_\gamma \xi$, with $c \in \mathcal C$ and $\gamma \in \Gamma$. Note that the set $\mathcal D$  is dense because $\xi$ is cyclic for $\pi({\mathcal C}\rtimes \Gamma)$. 

Now, it can be verified by a simple calculation that 
\bean
\lefteqn{\langle \pi_0(c_1) U_{\gamma_1} \xi, \pi_0(c_\alpha) \pi_0(c_2)U_{\gamma_2} \xi \rangle}\\
& =& \phi_{\gamma_1^{-1}\gamma_2}(\beta_{\gamma_1^{-1}}(c_1^* c_\alpha c_2)) \\
 & \rightarrow & \phi_{\gamma_1^{-1}\gamma_2}(\beta_{\gamma_1^{-1}}(c_1^* c c_2))\\
 &= & \langle \pi_0(c_1) U_{\gamma_1} \xi, \pi_0(c) \pi_0(c_2)U_{\gamma_2} \xi \rangle, \eean by the assumption of weak continuity of $\phi_\gamma$ for every $\gamma$. This completes the proof. 
\end{proof}

For $k \equiv (k_1,k_2) \in \MBBZ^2$, let    $\psi_{T}^k$ be the bounded linear functional  on the commutative $C^*$-algebra $C^*(U_l,V_r) \equiv {\C}(\MBBT^2)$ defined by
$$ \psi_{T}^k (U_l^m V_r^n)=\psi_{T} \left((U^m_l \otimes V_r^n)(V_l^{-1} \otimes V_r)^{k_1}(U_l^{-1} \otimes U_r)^{k_2}\right),~\forall m,n \in \MBBZ.$$
Let us denote by $\mu_T^k$ the complex measure on $\MBBT^2$ corresponding to $\psi_T^k$, i.e. $\psi_T^k(U_l^mV_r^n)=\int z_1^m z_2^n d\mu_T^k(z).$ From Proposition \ref{extend}, we conclude the following.

\bthm 
\label{ext}The CP map  $T$ is extendible  if and only if for each $k\in \MBBZ^2$ the measure $\mu_{T}^k$  is absolutly continuous with respect to Lebesgue measure on $\MBBT^2.$
\ethm

We end this subsection by mentioning a few examples of extendible maps. The verification of extendibility, using Theorem \ref{ext}, is quite straightforward, and we omit these calculations.    

\bxmpl   $T(a)=R^{*}aR,$ with $R \in \mathcal A,$ and is given by
$R=\sum_{i,j\in \MBBZ} c_{i,j}U_l^iV_l^j$ such that  $(c_{i,j})\in
l^2(\MBBZ^2).$ Then $T$ is extendible.
\exmpl

\bxmpl
Suppose that for each $k \in \MBBZ$, there exists  some positive integer  $
N_k$ such that $
T(f(U_l))V_l^k \in {\rm Span}[U_l^m V_l^n :m \in \MBBZ,~ -N_k \leq n \leq N_k],\ \forall f\in
{\C}(\MBBT).$ Then $T$ is extendible.
\exmpl

 \subsection{Invariants for extendible CP maps}

Suppose that $T$ is an extendible CP map, and let us denote the normal extension of $\psi_T$ again by the same symbol.
So, ${\psi_{T}}$ is a  normal state on  ${\mathcal A }
\otimes  {\mathcal A^\prime}$, which is also a
type $II_1$ hyperfinite factor, with  the canonical trace $\tau^{(2)}(\cdot)=\langle \xi_0, \cdot \xi_0 \rangle$ (where $\xi_0=1 \otimes 1 \in {\mathcal H}_\tau  \otimes {\mathcal H}_\tau$). Recall also the notation $c^{(r)}$ for $c \in \A \otimes \A^\prime$, introduced earlier by us. 

Since $\tau^{(2)}$ is faithful normal on $\A \otimes \A^\prime$, any normal state is absolutely continuous with restpect to $\tau^{(2)}$. Thus, we can obtain  a `Radon-Nikodym derivative', which is  a positive (possibly unbounded) 
operator $D_T$ affiliated to the commutant of $\A \otimes \A^\prime$, i.e.  $\A^\prime \otimes \A$,  
such that $\xi_0 \in {\rm Dom}(D_T^{\frac{1}{2}}),$ and 
\[\psi_{T}(b)=\tau^{(2)}(D_T b)\equiv \langle D_T^{\frac{1}{2}} \xi_0, b D_T^{\frac{1}{2}}\xi_0 \rangle,\ \forall b\in \A \otimes \A^\prime.\]
Note that $D_T$ is actually an element in $L^1(\A^\prime \otimes \A, \tau^{(2)})$.

We now explain how $D_T$ can be used to get numerical invariants for the equivalence $\sim$. Suppose that $S \sim T$, and let $u \in \A$ be such that $S(a)=T(u^*au)$ for all $a \in \A$. By an easy calculation we see that $$ \psi_S(b)=\psi_T(\tilde{u}^*b \tilde{u}),~~b \in \A \otimes \A^\prime,$$ where $\tilde{u}:=u \otimes 1$. Using the facts that $D_T^{\frac{1}{2}}$ commutes with $\tilde{u}$, and that $\tilde{u} \xi_0=\tilde{u}=\tilde{u}^{(r)} \xi_0,$ (where $\tilde{u}^{(r)} \in \A^\prime \otimes \A$ denotes right multiplication by $\tilde{u}$), we get the following :
\bean 
\lefteqn{\psi_S(b)}\\
& =& \psi_T(\tilde{u}^*b \tilde{u})\\
&=& \langle \tilde{u} D_T^{\frac{1}{2}} \xi_0, b \tilde{u} D_T^{\frac{1}{2}} \xi_0 \rangle\\
&=&  \langle  D_T^{\frac{1}{2}} \tilde{u}^{(r)} \xi_0, b  D_T^{\frac{1}{2}} \tilde{u}^{(r)}\xi_0 \rangle \\
& =& \langle (\tilde{u}^{(r)})^* D_T^{\frac{1}{2}} \tilde{u}^{(r)} \xi_0, b (\tilde{u}^{(r)})^* D_T^{\frac{1}{2}}\tilde{u}^{(r)} \xi_0 \rangle.
\eean
From this, we conclude (using the uniqueness of the Radon-Nikodym derivative) that $$ D_S=(\tilde{u}^{(r)})^*D_T \tilde{u}^{(r)}.$$ So, in particular, $\tau^{(2)}(e^{-tD_S})=\tau^{(2)}(e^{-tD_T})$ for all $t \geq 0$. These numbers can be used as invariants for $\sim$. Another possibility is to use $\tau^{(2)}(P_T)$, where $P_T$ denotes the projection (in $\A^\prime \otimes \A$) onto the range of $D_T$. 

\brmrk
Since $D_T$ is affiliated to  $\A^\prime \otimes \A$ , which can be identified with the von Neumann crossed product $L^\infty(\MBBT^2) \rtimes \MBBZ^2,$ we can describe $D_T$ by a formal series (which can be made sense of in a suitable $L^1$-topology) of the form $\sum_{k \in \MBBZ^2} D_T^{(k)} \delta_k$, where $D^{(k)}_T$ is a measurable complex valued function on $\MBBT^2$. Let us give a formula for the functions $D^{(k)}_T$ for a class of CP maps.   
 Consider $T(a)=R^* aR$ with  $R=\sum_{(m_1,m_2)}c_{m_1,m_2} U_l^{m_1} V_l^{m_2}$, where the summation is over a finite set, i.e. $c_m \equiv c_{(m_1,m_2)}$ is zero for all but a finitely many $m$'s. Then $T$ is extendible and $D^{(k)}_T$'s are given by,  
$$D_{T}^{(k)}(z)=\sum_{n,m \in \MBBZ^2 } c_n \bar{c}_{m_1+n_1,m_2+n_2+k_1-k_2}\lambda^{N(k,n,m)} z_1^{m_1}z_2^{m_2},$$  
  with $N(k,n,m)=2k_2 m_2+k_1 n_1 +k_2n_2+k_1k_2-k_2^2$. This can be verified by a simple calculation.

In particular, for $R=U$, 
$D_{T}^{(k)}(z)=\lambda^{k_2^2-k_1k_2}  z_2^{k_2-k_1},$  and 
for  $R=U+V$, 
$D_{T}^{(k)}(z)=\lambda^{k_2^2-k_1k_2}  z_2^{k_2-k_1}\{ \lambda^{k_1}+\lambda^{k_2}+\lambda^{2k_2+k_1} z_1^{-1}z_2+\lambda^{k_2}z_1^{-1}z_2^{-1}\}.$
\ermrk

Let us conclude with a nice computable formula for $\tau^{(2)}(P_T)$ when $D_T$ is $L^2$.

\bpro If $D_T$ belongs to $L^2(\A^\prime \otimes \A, \tau^{(2)}), $ then we have :
\[{\tau^{(2)}}(P_T)=1-\inf_{b\in  {\A^{\rm fin}}}\  \sup_{c\in \A^{\rm fin}:\tau^{(2)}(c^*c)=1} |\tau^{(2)}( c^*)-\psi_T(c^*b)|^2.\]
\epro
\begin{proof}
Note that $D_T \in L^2$ implies that $\xi_0 \in {\rm Dom}(D_T)$, and thus, $\tau^{(2)}(b)=\langle \xi_0, b D_T \xi_0 \rangle$ for all $b \in \A \otimes \A^\prime$. Moreover, As $D_T $ is affiliated to $(\A \otimes \A^\prime)^{\prime}$, 
$ e^{-tD}$ leaves $\A^{\rm fin}\xi_0 \subseteq \mathcal H_\tau \otimes \mathcal H_\tau$ invariant. Hence $\A^{\rm fin} \xi_0$ is a core for $D_T$, hence  
$\{ D_T b \xi_0~:~b \in \A^{\rm fin} \}$ is dense in ${\rm Ran}~D$. Thus, we have, 
\bean \lefteqn {\tau^{(2)}(P_T)}\\
&=& \| P_T\ \xi_0 \|^2\\
&=& 1-dist(\xi_0,{\rm Ran}~D_T)^2\\
&=& 1-\inf_{b\in  {\A^{\rm fin}}}\| \xi_0-D_T\ b \xi_0 \|^2  \\
&=& 1-\inf_{b\in  {\A^{\rm fin}}}\  \sup_{c\in \A^{\rm fin}:\tau^{(2)}(c^*c)=1}
|\langle c,\xi_0-D_T  b \xi_0 \rangle|^2  \\
&=& 1-\inf_{b\in  {\A^{\rm fin}}}\  \sup_{c\in \A^{\rm fin}:\tau^{(2)}(c^*c)=1}
|\tau^{(2)}( c^*)-\tau^{(2)}(c^*\ D_T b)|^2  \\
&&=1-\inf_{b\in  {\A^{\rm fin}}}\  \sup_{c\in \A^{\rm fin}:\tau^{(2)}(c^*c)=1} |\tau^{(2)}( c^*)-\psi_{T}(c^*b)|^2 \eean
\end{proof}

However, in this paper, we do not study these numerical invariants in more details, and do not investigate  whether one or more of such invariants characterize the equivalence $\sim$.  We leave these topics for future research.  
\vspace{3mm}\\
\noindent
{\bf Appendix}\\
Here we give an example of a commuting pair of normal representations $\pi$ and $\rho$ of  a von Neumann algebra $\mathcal A \subseteq {\mathcal B}({\mathcal H}_0)$ and its commutant respectively, such that the representation $\Phi$ of $\A \otimes_{\rm alg} \A^\prime$ defined by $\Phi(a \otimes a^\prime)=\pi(a) \rho(a^\prime)$ does not have a normal extension to  the von Neumann algebra $\A \otimes \A^\prime$. To this end, let us take ${\mathcal H}_0=L^2(\MBBT)$ (with Lebesgue measure), and $\A=\A^\prime=L^\infty(\MBBT)$. Let $\mathcal H=L^2(\MBBR^2, \mu)$, where $$ \mu(E):= l(E \bigcap D),$$ $D := \bigcup_{n \in \MBBZ} \{ (t,t+n)~:~t \in \MBBR \}$, and $l$ denotes the Lebesgue measure on $\MBBR$. Define $\pi$ and $\rho$ from $L^\infty(\MBBT)$ to $\mathcal B(\mathcal H)$ by setting : $$ (\pi(\phi)f)(s,t)=\phi(\lambda^s)f(s,t),~~(\rho(\phi)f)(s,t)=\phi(\lambda^t)f(s,t),$$ for $f \in L^2(\MBBR^2,\mu)$. We claim that $\Phi : L^\infty(\MBBT) \otimes_{\rm alg} L^\infty(\MBBT) \rightarrow {\mathcal B}(\mathcal H)$ given by $\Phi(\phi \otimes \psi):=\pi(\phi)\rho(\psi)$, does not admit a normal extension to $L^\infty(\MBBT^2)$. To see this, note that the set $\Gamma :=\{ (x,(\lambda^n x){~\rm mod}~1)~:~n \in \MBBZ \}$ has zero Lebesgue measure as a subset of $\MBBT^2$. Now, assuming, if possible, that $\Phi$ admits a normal extension on $L^\infty(\MBBT^2)$, choose a sunction $F \in {\mathcal C}(\MBBT^2)$ and define $F_0 $ in $L^\infty(\MBBT^2)$ by setting $F_0(z)=0$ for $z \in \Gamma$, and $F_0(z)=F(z)$ for all $z$ in the complement of $\Gamma$. Thus, $F=F_0$ a.e. (Lebesgue), hence $\Phi(F)=\Phi(F_0)$. However, \bean
\lefteqn{\| \Phi(F_0) f \|^2}\\
&=& \int_{\MBBR^2} | F_0(\lambda^s, \lambda^t) f(s,t)|^2 d \mu(s,t)\\
&=& \sum_{n \in \MBBZ} \int_{\MBBR} |F_0(\lambda^s,\lambda^n \lambda^s) f(s,s+n) |^2 ds \\
&=& 0,\\
\eean
for every $f \in L^2(\mu)$. Thus, $\Phi(F)=0$ for all $F \in {\mathcal C}(\MBBT^2)$, contradicting the assumption that $\Phi$ extends normally.

\end{document}